\documentclass[11pt,twoside]{amsart}
\usepackage{amsxtra}
\usepackage{color}
\usepackage{amsopn}
\usepackage{amsmath,amsthm,amssymb}

\title{Half-Flat Structures Inducing Einstein Metrics on Homogeneous Spaces}
\author{Alberto Raffero}
\address{Dipartimento di Matematica ``G. Peano", Universit\`a di Torino, via Carlo Alberto 10, 10123 Torino, Italy}
\email{alberto.raffero@unito.it}

\theoremstyle{remark}
\newtheorem{remark}{Remark}[section]

\theoremstyle{definition}
\newtheorem{defin}[remark]{Definition}
\newtheorem{ex}[remark]{Example}

\theoremstyle{plain}
\newtheorem{teo}[remark]{Theorem}
\newtheorem{prop}[remark]{Proposition}
\newtheorem{corol}[remark]{Corollary}
\newtheorem{lemma}[remark]{Lemma}

\newcommand{\beq}{\begin{equation}}
\newcommand{\eeq}{\end{equation}}
\newcommand{\bqn}{\begin{eqnarray}}
\newcommand{\eqn}{\end{eqnarray}}
\newcommand{\bqne}{\begin{eqnarray*}}
\newcommand{\eqne}{\end{eqnarray*}}

\newcommand{\R}{{\mathbb R}}

\newcommand{\C}{{\mathbb C}}

\newcommand{\W}{\wedge}

\newcommand{\psip}{\psi_+}
\newcommand{\psim}{\psi_-}
\newcommand{\ff}{{\rm f}}
\newcommand{\SU}{{\rm SU}}

\begin{document}
\maketitle
\begin{abstract}
In this paper, we consider half-flat $\SU(3)$-structures and the subclasses of coupled and double structures. In the general case we show that the intrinsic torsion form $w_1^-$ is constant in each of the two subclasses. We then consider the problem of finding half-flat structures inducing Einstein metrics on homogeneous spaces. We give an example of a left invariant half-flat (non coupled and non double) structure inducing an Einstein metric on $S^3\times S^3$ and we show there does not exist any left invariant coupled structure inducing an ${\rm Ad}(S^1)$-invariant Einstein metric on it. Finally, we show that there are no coupled structures inducing the Einstein metric on Einstein solvmanifolds and on homogeneous Einstein manifolds of nonpositive sectional curvature.
\end{abstract}

\section{Introduction}
An $\SU(3)$-structure on a 6-dimensional smooth manifold $N$ is the data of a Riemannian metric $h$, an orthogonal almost complex structure $J,$ a 2-form $\omega$ related with $h$ and $J$ via the identity $\omega(\cdot,\cdot)=h(J\cdot,\cdot)$ and a $(3,0)$-form $\Psi$ of nonzero constant length. Using the results on stable forms (\cite{Hit}, \cite{Rei}) it can be shown that such a structure actually depends only on $\omega$ and $\psip:=\Re(\Psi)$.

The intrinsic torsion of an $\SU(3)$-structure lies in a 42-dimensional space whose decomposition into $\SU(3)$-irreducible summands allows to divide the $\SU(3)$-structures in classes, which can be described by the exterior derivatives of $\omega, \psip$ and $\psim:=\Im(\Psi)=J^*\psip$ as shown in \cite{CS} or in terms of a characterizing spinor and the spinorial field equations it satisfies as recently shown in \cite{ACFH}.

 When the intrinsic torsion vanishes, i.e.~when $\omega,\psip$ and $\psim$ are all closed, the manifold has holonomy contained in $\SU(3)$ and the metric $h$ is Ricci-flat.

When the forms $\psip$ and $\omega^2 = \omega\W\omega$ are both closed, the torsion lies in the 21-dimensional space $\mathcal{W}_1^-\oplus\mathcal{W}_2^-\oplus\mathcal{W}_3$ and the $
\SU(3)$-structure is called {\it half-flat}. Half-flat structures are the initial values for the Hitchin flow equations and are used as starting point to construct 7-dimensional manifolds with holonomy 
contained in $G_2$ (see for example \cite{Bry}, \cite{CS}, \cite{CLSS}, \cite{Hit}).

There are three classes of $\SU(3)$-structures contained in the half-flat one which have been frequently considered in literature: the {\it nearly K\"ahler} ($\mathcal{W}_1^-$), the {\it double} or 
{\it co-coupled} ($\mathcal{W}_1^-\oplus\mathcal{W}_3$) and the {\it coupled} ($\mathcal{W}_1^-\oplus\mathcal{W}_2^-$) (see \cite{Sal}). 
It is well known (\cite{Gr4},\cite{Ma}) that the metric $h$ induced by a nearly K\"ahler structure is always Einstein, i.e.~the Ricci tensor of $h$ satisfies the identity
$${\rm Ric}(h)  = \mu h$$
for some real number $\mu$. Moreover there exist some examples of half-flat structures with torsion class $\mathcal{W}_1^-\oplus\mathcal{W}_2^-\oplus\mathcal{W}_3$ and $\mathcal{W}_1^-\oplus\mathcal{W}_3$ inducing an Einstein metric, but up to now it seems the coupled case has not been studied. 

The almost Hermitian structure $(h, J, \omega)$ underlying a coupled $\SU(3)$-structure is quasi K\"ahler and depends only on $\omega$. Thus the properties of quasi K\"ahler 
manifolds, e.g.~the curvature identities shown by Gray in \cite{Gr3}, are still valid on manifolds endowed with a coupled structure.
The class of quasi K\"ahler structures contains the almost K\"ahler one. It was conjectured by Goldberg in \cite{Go} that any compact almost K\"ahler manifold whose associated metric is Einstein is actually K\"ahler. In \cite{Se}, Sekigawa showed that this conjecture is true when the scalar curvature is non negative. Moreover, there exists a non compact example of an Einstein almost K\"ahler manifold with negative scalar curvature (\cite{ADM}), which is the unique example of 6-dimensional Einstein almost K\"ahler (non K\"ahler) solvmanifold by the results 
contained in \cite{HOS} and \cite{FFM}.

In \cite{FR}, invariant coupled structures on 6-dimensional nilpotent Lie groups were classified and it was shown that there is only one case in which the coupled structure induces a Ricci soliton metric. Moreover, this coupled structure was used to construct a locally conformal calibrated $G_2$-structure inducing an Einstein metric on the rank one solvable extension of the nilpotent Lie algebra.

More in general, it is not difficult to show that on the cylinder and on the cone over a 6-manifold admitting a coupled structure there exists a locally conformal calibrated $G_2$-structure induced by the coupled one. On the other hand, it is possible to show that a parallel $G_2$-structure on a 7-dimensional manifold induces a coupled structure $(h, J, \omega, \Psi)$ on the oriented hypersurfaces having $J$-invariant second fundamental form (\cite{Cab}).

One of the simplest cases which can be considered when one is looking for examples of special geometric structures with (or without) torsion is the one of left invariant structures on homogeneous spaces, since in this case the starting analytic problem on the manifold (e.g.~the problem of solving the PDEs deriving from the definition of half-flat structure) can often be reduced to an algebraic problem on the tangent space to a point. Following this idea, this paper focuses the attention to the case of left invariant half-flat structures on 6-dimensional homogeneous  manifolds.

In \cite{Sch}, Schulte-Hengesbach considered half-flat structures on Lie groups and, in particular, he described left invariant half-flat structures on $S^3\times S^3$ and gave two examples of half-flat structures inducing an Einstein metric on it. One of these examples consists of a double structure and the other is the unique (up to homotheties and sign) left invariant nearly K\"ahler structure existing on this manifold, as showed by Butruille in \cite{Bu}. 
Moreover, the Einstein metrics are the only two currently known examples of left invariant Einstein metrics on $S^3\times S^3$ and can be characterized as the only ${\rm Ad}(S^1)$-invariant Einstein metrics existing on it up to isometries and homotheties.
It is then natural to ask whether there exist left invariant half-flat structures (neither coupled nor double) and left invariant coupled structures on $S^3 \times S^3$ inducing any of these metrics. We  show that in the first case the answer is positive by giving an explicit example, while for the coupled case we prove that the answer is negative. 

Left invariant half-flat structures on $S^3\times S^3$ were also studied by Madsen and Salamon in \cite{MS}, where they described them using the representation theory of ${\rm SO}(4)$ and matrix algebra and showed that the moduli space they define is essentially a finite dimensional symplectic quotient. In that paper they also considered the subclasses of coupled, co-coupled (double) and nearly K\"ahler structures, in particular they gave some examples of double structures, gave another proof of Butruille's result in their setting and constructed a 1-parameter family of double structures which is a solution for the Hitchin flow.

Conversely to the case of the compact manifold $S^3\times S^3$, in the case of noncompact homogeneous manifolds it can be shown using the results of Heber (\cite{He}) and Lauret (\cite{Lau}) that every Einstein solvmanifold has a unique left invariant Einstein metric up to isometries and homotheties. Einstein solvmanifolds constitute the unique example of noncompact homogeneous Einstein manifolds known up to now and it has been conjectured by Alekseevskii that they might be the only case that can occur (\cite[7.57]{Bes}). 

We then focus on 6-dimensional Einstein solvmanifolds, which were classified by Nikitenko and Nikonorov in \cite{NN}, and we look for left invariant half-flat structures defined on them inducing the Einstein metric. In particular, the existence of a coupled structure inducing the Einstein metric would provide an example of a solvmanifold endowed with 
a left invariant quasi K\"ahler Einstein structure. However, we prove that in this case there are no coupled structures satisfying the property we are looking for.
Moreover, using another result contained in \cite{NN}, we are able to conclude that this happens for 6-dimensional homogeneous Einstein manifolds of nonpositive sectional curvature too. This gives a constraint to the existence of coupled Einstein structures on 6-dimensional homogeneous manifolds and leads to ask whether the non positivity of the sectional curvature gives a constraint in a more general setting. 

This paper is organized as follows: in section \ref{sechf} we recall the definition of $\SU(3)$-structures and we focus on the half-flat class and its subclasses, 
in section \ref{s3s3} we consider the case of left-invariant half-flat structures on $S^3\times S^3$ and in section \ref{secsol}, after recalling some properties of 6-dimensional Einstein solvmanifolds, we study the case of left invariant half-flat structures inducing Einstein metrics on them.

All the algebraic computations in section \ref{s3s3} and \ref{secsol} have been done with the aid of the software Maple.

\medskip
\noindent {\em{Acknowledgements}}.  The author would like to thank Anna Fino for suggesting the problem and for the fruitful conversations on the work and Thomas Madsen and Paolo Lella for the useful comments.

\section{Preliminaries on $\SU(3)$-structures}\label{sechf}
Let $N$ be a 6-dimensional smooth manifold. It admits an $\SU(3)$-structure if the structure group of its frame bundle can be reduced from ${\rm GL}(6,\R)$ to $\SU(3)$. 
The existence of an $\SU(3)$-structure is equivalent to the existence of an almost Hermitian structure $(h, J, \omega)$ and a $(3,0)$-form $\Psi$ of nonzero constant length such that 
the Riemannian metric $h$, the orthogonal almost complex structure $J$ and the 2-form $\omega$ are related via the identity 
$$\omega(\cdot,\cdot)=h(J\cdot,\cdot),$$
the forms $\omega$, $\psip:=\Re(\Psi)$ and $\psim:= \Im(\Psi) = J^*\psip$ are {\bf compatible} in the sense that
$$\omega\W\psi_\pm=0$$
and satisfy the {\bf normalization condition}
$$\psip\W\psim=\frac23\omega^3 = 4dV_h,$$
where ${\rm d}V_h$ is the Riemannian volume form of $h$.

Using the results on stable forms contained in the works \cite{Hit}, \cite{Rei}, it is possible to show that an $\SU(3)$-structure actually depends only on $\omega$ and $\psip$ and thus it is possible to give another characterization for $\SU(3)$-structures, which is the one we use in this work. To describe this characterization let us consider a 6-dimensional real vector space $V,$ we say that a differential $k$-form on $V$ is {\bf stable} if its orbit under the natural action of ${\rm GL}(V)$ on $\Lambda^k(V^*)$ is open. One can show that a 2-form $\sigma$ on $V$ is stable if and only if $\sigma^3\neq0$, that is if and only if it is non degenerate. Moreover, if we 
denote by $A:\Lambda^5(V^*) \rightarrow V\otimes\Lambda^6(V^*)$ the canonical isomorphism given by $A(\gamma) = v \otimes \Omega$, where $i_v\Omega = \gamma$, 
and define for a fixed $\rho \in\Lambda^3(V^*)$
$$K_\rho : V \rightarrow V\otimes\Lambda^6(V^*),\ \  K_\rho(v) = A((i_v \rho)\W\rho)$$ 
and 
$$\lambda : \Lambda^3(V^*) \rightarrow (\Lambda^6(V^*))^{\otimes2},\ \  \lambda(\rho) = \frac16{\rm tr}K^2_\rho,$$
we have that a 3-form $\rho$ is stable if and only if $\lambda(\rho)\neq0$. Whenever this happens, one can choose the orientation of $V$ for which $\omega^3$ is positively oriented and define a volume form by $\sqrt{|\lambda(\rho)|} \in \Lambda^6(V^*)$ and an endomorphism 
$$J_\rho = -\frac{1}{\sqrt{|\lambda(\rho)|}}K_\rho,$$ 
which is an almost complex structure when $\lambda(\rho)<0$. 

The existence of an $\SU(3)$-structure on $N$ is equivalent to the existence of a pair of differential forms $(\omega,\psip)\in\Lambda^2(N)\times\Lambda^3(N)$ such that for each $p\in N$ the forms $\omega(p)$ 
and $\psip(p)$ on $T_pN$ are stable with $\lambda(\psip(p))<0$, compatible, satisfy the normalization condition and define a Riemannian metric 
$h_p(\cdot,\cdot)=\omega_p(\cdot,J\cdot)$, where $J = J_{\psip(p)}$ is the almost complex structure induced by $\psip(p)$.
The Riemannian metric $h$ can also be described in terms of $\omega$ and $\psip$ as 
$$h(X,Y)\omega^3 = -3(i_X\omega)\W(i_Y\psip)\W\psip,$$ 
for any $X,Y\in\mathfrak{X}(N)$.

The intrinsic torsion of an $\SU(3)$-structure lies in a 42-dimensional space 
$$\mathcal{W}_1^+\oplus\mathcal{W}_1^-\oplus\mathcal{W}_2^+\oplus\mathcal{W}_2^-\oplus\mathcal{W}_3\oplus\mathcal{W}_4\oplus\mathcal{W}_5$$
given by the sum of irreducible $\SU(3)$-modules and completely determined by $d\omega, d\psip$ and $d\psim$ (see \cite{CS}).

In this paper we are mainly interested in $\SU(3)$-structures having both $\omega^2$ and $\psip$ closed, known as {\bf half-flat} $\SU(3)$-structures in literature. 
In this case the intrinsic torsion lies in the space $\mathcal{W}_1^-\oplus\mathcal{W}_2^-\oplus\mathcal{W}_3$ and in terms of the exterior derivatives of $\omega, \psip$ and $\psim$ this reads:
\begin{equation}
\begin{array}{lcl}\label{hf}
d\omega &=& -\frac32w_1^-\psip  + w_3,\\
d\psip &=&0,\\
d\psim &=& w_1^-\omega^2 - w_2^-\W\omega,
\end{array}
\end{equation}
where $w_1^-\in C^\infty(N)\cong \mathcal{W}_1^-, w_2^-\in\Lambda^{1,1}_0(N)\cong\mathcal{W}_2^-, w_3 \in \Lambda^{2,1}_{0}(N)\cong\mathcal{W}_3$ are the (non vanishing) 
{\bf intrinsic torsion forms}.

There are three interesting families of $\SU(3)$-structures with nonzero torsion contained in the family of half-flat ones, we recall here the definitions.
\begin{defin}
A half-flat $\SU(3)$-structure is said to be {\bf nearly K\"ahler} if $\nabla J$ is skew-symmetric, that is $\nabla_X(J)(X) = 0$ for every $X\in\mathfrak{X}(N)$, {\bf coupled} if $d\omega  = c\psip$ for some non vanishing $c\in C^\infty(N)$ and {\bf double} (or {\bf co-coupled}) if $d\psim  = k\omega^2$ for some non vanishing $k\in C^\infty(N)$.
\end{defin}
It can be shown that the intrinsic torsion of a nearly K\"ahler structure lies in $\mathcal{W}_1^-$, while it follows easily from the definition that the intrinsic torsion of a coupled  
lies in $\mathcal{W}_1^-\oplus\mathcal{W}_2^-$ and the intrinsic torsion of a double lies in $\mathcal{W}_1^-\oplus\mathcal{W}_3$. 
In each case the expressions of $d\omega, d\psip$ and $d\psim$ can be obtained from \eqref{hf} having in mind the torsion class to which each family belongs. Looking at these, it 
is easy to see that coupled and double structures can be thought as a generalization of the nearly K\"ahler.
For each class it is also possible to write the Ricci tensor and the scalar curvature of $h$ in terms of the non vanishing intrinsic torsion forms using the results of \cite{BV}. 
One can then recover that the metric induced by a nearly K\"ahler structure is Einstein and can observe that in the general case this is not true anymore for the metric induced by a coupled or a double.

It is well known that in the nearly K\"ahler case the only non vanishing torsion form $w_1^-$ is constant, using the fact that the coupled and the double structures are in particular half-flat, we can prove that the same is true in these two cases. 
\begin{lemma}
Let $N$ be a 6-dimensional connected smooth manifold endowed with an $\SU(3)$-structure $(\omega,\psip)$. 
If $(\omega,\psip)$ is coupled, then there exists a nonzero real constant $c$ such that 
$$d\omega = c\psip,$$
if $(\omega,\psip)$ is double, then there is a nonzero real constant $k$ such that
$$d\psim = k\omega^2.$$
\end{lemma} 
\proof
If $(\omega,\psip)$ is coupled, using the notations of \eqref{hf} we know that
$$d\omega =-\frac32w_1^-\psip,$$
where $w_1^-$ is a smooth nonzero function.
Taking the exterior derivative of $d\omega$ we obtain
\begin{eqnarray}
0 &=& d(w_1^-\psip)\nonumber\\
   &=& dw_1^-\W\psip + w_1^- d\psip,\nonumber
\end{eqnarray}
so
$$d\psip = -\frac{1}{w_1^-}dw_1^-\W\psip.$$
Now, $d(\omega^2) = 2d\omega\W\omega = 0$ since $\omega$ and $\psip$ are compatible, thus the considered class of $\SU(3)$-structures is contained in the half-flat one if and only if
$$dw_1^-=0,$$
that is $w_1^-$ is a nonzero real constant on $N$ connected. 

In the double case we can argue in a similar way: starting from 
$$d\psim = w_1^-\omega^2,$$
we take the exterior derivative of both sides obtaining
$$0 = dw_1^-\W\omega^2 + w_1^-d\omega^2$$
and conclude observing that $d\omega^2 = 2d\omega\W\omega = 0$ since $\omega\W\psip = 0 = \omega\W w_3$ and that wedging 1-forms by $\omega^2$ is an isomorphism.
\endproof

\begin{remark}
The defining conditions of a half-flat structure $d\psip=0$ and $d\omega^2=0$ are obviously satisfied by the subclasses of coupled, double and nearly K\"ahler, thus  
to make distinction between these classes is necessary to look at the non vanishing components of the intrinsic torsion (for example by computing $d\omega$ and $d\psim$). 
Moreover, we could sometimes emphasize the fact that we are considering half-flat structures with torsion in $\mathcal{W}_1^-\oplus\mathcal{W}_2^-\oplus\mathcal{W}_3$ by saying that 
the structure is half-flat non coupled and non double.
\end{remark}

Recall that an almost Hermitian structure $(h,J,\omega)$ on a $2n$-dimensional manifold $M$ is said to be {\bf quasi K\"ahler} if $\bar\partial\omega = (d\omega)^{1,2} = 0$ or, equivalently, if for any $X,Y\in\mathfrak{X}(M)$ it holds $\nabla_X(J)(Y)+\nabla_{JX}(J)(JY)=0,$ where $\nabla$ is the Levi Civita connection of $h$. 
Properties of quasi K\"ahler manifolds have been studied in literature (see for example \cite{Gr1},\cite{Gr3},\cite{VH}) and in dimension 6 are inherited by  
manifolds endowed with a coupled structure. Indeed we have
\begin{prop}
Let $(\omega,\psip)$ be a coupled $\SU(3)$-structure on a 6-dimensional connected smooth manifold $N,$ then the underlying almost Hermitian structure $(h,J,\omega)$ is quasi K\"ahler.
Moreover it depends only on $\omega$.
\end{prop}
\proof
As we already observed, we have that $d\omega = c\psip$ for some $c\in\R-\{0\}$. Thus $d\omega$ is of type $(3,0)+(0,3)$ and $(d\omega)^{1,2}=0$. 
Moreover, since $J$ depends only on $\psip$, which is proportional to $d\omega$, the almost Hermitian structure $(h,J,\omega)$ depends only on $\omega$.
\endproof

\section{Einstein half-flat structures on $S^3\times S^3$}\label{s3s3}
In this section we focus on the homogeneous manifold $S^3\times S^3$ and we consider left invariant half-flat structures on it inducing Einstein metrics, looking for what happens in the subclasses.
 
As a Lie group the manifold is $\SU(2)\times \SU(2)$ and we can describe the left invariant tensors only working on the Lie algebra $\mathfrak{su}(2)\oplus\mathfrak{su}(2)$. In particular, 
a left invariant metric on $S^3\times S^3$ can be identified with an inner product on  $\mathfrak{su}(2)\oplus\mathfrak{su}(2)$
and a left invariant $\SU(3)$-structure $(\omega,\psip)$ on $S^3\times S^3$ can be identified with a 2-form $\omega$ and a 3-form $\psip$ defined on $\mathfrak{su}(2)\oplus\mathfrak{su}(2)$  satisfying the defining properties of such a structure.

Let $(e_1, e_2, e_3)$ denote the standard basis for the first copy of $\mathfrak{su}(2)$, $(e_4, e_5, e_6)$ denote it for the second one and let $(e^1,e^2,e^3)$ and $(e^4,e^5,e^6)$ denote their dual bases. Then we have the following structure equations for the Lie algebra $\mathfrak{su}(2)\oplus\mathfrak{su}(2)$:
$$\renewcommand\arraystretch{1.4}
\begin{array}{c}
de^1 = e^{23}, de^2 = e^{31}, de^3 = e^{12},\\
de^4 = e^{56}, de^5 = e^{64}, de^6 = e^{45},
\end{array}
$$
where we are using the shortening $e^{ijk\cdots}$ for the wedge product $e^i\W e^j \W e^k \W\cdots$.
 
The problem of classifying left invariant Einstein metrics on $S^3\times S^3$ is still open and only two examples are known. These can be characterized as follows:
\begin{teo}[\cite{NR}]\label{nikeinst}
Let $h$ be a left invariant Einstein metric on the Lie group $\SU(2)\times \SU(2)$ which is ${\rm Ad}(S^1)$-invariant for some embedding $S^1\subset \SU(2)\times \SU(2)$,  
then $h$ is isometric up to homotheties either to the standard metric or to Jensen's metric. 
\end{teo}

With respect to the basis $(e_1,e_2,e_3,e_4,e_5,e_6)$ and up to scalar multiples, the matrix associated to the standard metric is the identity matrix and the one associated to the Jensen metric is the following
\begin{equation}\label{Jensen}
\left[ \begin {array}{cccccc} \frac{2\sqrt{3}}{3}\,&0&0&-\frac{\sqrt{3}}{3}\,&0&0
\\ \noalign{\medskip}0&\frac{2\sqrt{3}}{3}\,&0&0&-\frac{\sqrt{3}}{3}\,&0
\\ \noalign{\medskip}0&0&\frac{2\sqrt{3}}{3}\,&0&0&-\frac{\sqrt{3}}{3}\,
\\ \noalign{\medskip}-\frac{\sqrt{3}}{3}\,&0&0&\frac{2\sqrt{3}}{3}\,&0&0
\\ \noalign{\medskip}0&-\frac{\sqrt{3}}{3}\,&0&0&\frac{2\sqrt{3}}{3}\,&0
\\ \noalign{\medskip}0&0&-\frac{\sqrt{3}}{3}\,&0&0&\frac{2\sqrt{3}}{3}\,
\end {array} \right] .
\end{equation}

In \cite{Sch}, Schulte-Hengesbach gave an example of a left invariant half-flat structure on $S^3\times S^3$ inducing the standard metric and of one inducing the Jensen metric, we describe them in the following examples.
\begin{ex}[\cite{Sch}]
The pair of stable forms
$$\renewcommand\arraystretch{1.4}
\begin{array}{lcl}
\omega & = & -e^{14}-e^{25}-e^{36},\\
\psip & = & \frac{1}{\sqrt{2}}\left(e^{123}-e^{156}+e^{246}-e^{345}+e^{126}-e^{135}+e^{234}-e^{456}\right),
\end{array}
$$
is compatible, normalized and induces the standard metric, thus it defines an $\SU(3)$-structure on $\mathfrak{su}(2)\oplus\mathfrak{su}(2)$. 
Moreover $d\psip=0$, $d\omega^2=0$, $d\psim = \frac{1}{\sqrt{2}}\omega^2$ and $d\omega$ is not proportional to $\psip$, i.e.~it is a double structure.
\end{ex} 
\begin{ex}[\cite{Sch}]
The following pair of stable, compatible, normalized forms
$$\renewcommand\arraystretch{1.4}
\begin{array}{lcl}
\omega & = & -\frac{\sqrt{3}}{18}\left(e^{14}+e^{25}+e^{36}\right),\\
\psip &=& \frac{\sqrt{3}}{54}\left(-e^{234}+e^{156}+e^{135}-e^{246}-e^{126}+e^{345}\right),
\end{array}
$$
induces the Jensen metric, therefore it defines an $\SU(3)$-structure on $\mathfrak{su}(2)\oplus\mathfrak{su}(2)$ which is nearly K\"ahler since $d\omega = 3\psip$ and $d\psim=-2\omega^2$.
\end{ex}
\begin{remark}
The signs of the differential forms in the previous examples are different from those in the original examples of \cite{Sch}, 
this is due to the fact that in this paper we are using a convention in the definition of an $\SU(3)$-structure which is slightly different from the one used by the other author in his works. 
\end{remark}
We can also give a new example of a left invariant half-flat structure (non coupled and non double) inducing the Jensen metric:
\begin{ex}
The pair
$$\renewcommand\arraystretch{1.4}
\begin{array}{lcl}
\omega &=& \frac{\sqrt[3]{4}\sqrt[6]{3}}{2}\left(-e^{14}+e^{25}+e^{36}\right),\\
\psip &=& e^{123}+e^{135}-e^{246}-e^{126}+e^{345}-e^{456},
\end{array}
$$
defines an $\SU(3)$-structure on $\mathfrak{su}(2)\oplus\mathfrak{su}(2)$ and induces a metric which is (proportional to) the Jensen metric. Moreover, it can be checked that this $\SU(3)$-structure is half-flat since both $\psip$ and $\omega^2$ are closed and it is neither coupled nor double since $d\omega$ is not proportional to $\psip$ and $d\psim$ is not proportional to $\omega^2$.
\end{ex}

Summarizing, on $S^3\times S^3$ there exist left invariant half-flat and nearly K\"ahler structures inducing the Jensen metric and left invariant double structures inducing the standard metric. 
We prove now that it is not possible to find left invariant coupled structures on $S^3\times S^3$ inducing either the standard or the Jensen metric. In the proof we use some classical 
properties of algebraic varieties, the reader can find them for example in \cite{CLO}. 
\begin{teo}\label{nocpds3}
$S^3\times S^3$ does not admit left invariant coupled $\SU(3)$-structures $(\omega, \psip)$ inducing an ${\rm Ad}(S^1)$-invariant Einstein metric.
\end{teo}
\proof
Let us consider a left invariant coupled structure $(\omega,\psip)$ on $S^3\times S^3$ which we identify with a 2-form $\omega$ and a 3-form $\psip$ defined on $\mathfrak{su}(2)\oplus\mathfrak{su}(2)$ and such that $d\omega = \frac1c\psip$, $c\in\R-\{0\}$.
Since $\omega^2$ is closed, it follows that $\omega\in \mathfrak{su}^*(2)\otimes\mathfrak{su}^*(2)$ (cf.~\cite{Sch}, Lemma 1.1 p. 81), thus
$$\omega = a_{14}e^{14} + a_{15}e^{15} + a_{16}e^{16} + a_{24}e^{24} + a_{25}e^{25} + a_{26}e^{26} + a_{34}e^{34} + a_{35}e^{35} + a_{36}e^{36},$$
where $a_{ij}$ are real coefficients.
Imposing the coupled condition $\psip = cd\omega, c\neq 0$, we have that 
\begin{eqnarray}
\psip & = & c(a_{1 4}e^{234}-a_{1 4}e^{156}+a_{1 5}e^{235}+a_{1 5}e^{146}+a_{1 6}e^{236}-a_{1 6}e^{145}-a_{2 4}e^{134}-a_{2 4}e^{256}\nonumber\\
          &  & -a_{2 5}e^{135}+a_{2 5}e^{246}-a_{2 6}e^{136}-a_{2 6}e^{245}+a_{3 4}e^{124}-a_{3 4}e^{356}+a_{3 5}e^{125}+a_{3 5}e^{346}\nonumber\\
          &  & +a_{3 6}e^{126}-a_{3 6}e^{345})\nonumber
\end{eqnarray}
and from the closedness of $\omega^2$ we know that the compatibility condition $\omega\W\psip=0$ holds.
It is now possible to compute $\lambda = \lambda(\psip)$, which turns out to be a homogeneous polynomial of degree 4 in the coefficients $a_{ij}$, 
the almost complex structure $J = J_{\psip}$ and $h(\cdot,\cdot) = \omega(\cdot,J\cdot)$. 
With respect to the basis $(e_1,\ldots,e_6)$ the matrix $H$ associated to $h$ is symmetric. Moreover, up to a global sign depending on whether the considered basis is positively oriented or not and not affecting the computations afterwards, the nonzero entries are the following:
$$
\begin{array}{rcl}
H_{i,i}     & = & \frac{-2c^2}{\sqrt{-\lambda}}(a_{14}a_{25}a_{36}-a_{14}a_{26}a_{35}-a_{15}a_{24}a_{36}+a_{15}a_{26}a_{34}+a_{16}a_{24}a_{35}-a_{16}a_{25}a_{34}),\\
               & &i = 1,\ldots,6\\
H_{1,4}  & = & \frac{-c^2}{\sqrt{-\lambda}}(a_{14}^3+a_{14}a_{15}^2+a_{14}a_{16}^2+a_{14}a_{24}^2-a_{14}a_{25}^2-a_{14}a_{26}^2+a_{14}a_{34}^2-a_{14}a_{35}^2-a_{14}a_{36}^2\\
               & &+2a_{15}a_{24}a_{25}+2a_{15}a_{34}a_{35}+2a_{16}a_{24}a_{26}+2a_{16}a_{34}a_{36}),\\
H_{1,5}  & = & \frac{-c^2}{\sqrt{-\lambda}}(a_{14}^2a_{15}+2a_{14}a_{24}a_{25}+2a_{14}a_{34}a_{35}+a_{15}^3+a_{15}a_{16}^2-a_{15}a_{24}^2+a_{15}a_{25}^2-a_{15}a_{26}^2\\
               & &-a_{15}a_{34}^2+a_{15}a_{35}^2-a_{15}a_{36}^2+2a_{16}a_{25}a_{26}+2a_{16}a_{35}a_{36}),\\
H_{1,6}  & = & \frac{-c^2}{\sqrt{-\lambda}}(a_{14}^2a_{16}+2a_{14}a_{24}a_{26}+2a_{14}a_{34}a_{36}+a_{15}^2a_{16}+2a_{15}a_{25}a_{26}+2a_{15}a_{35}a_{36}+a_{16}^3\\
               & &-a_{16}a_{24}^2-a_{16}a_{25}^2+a_{16}a_{26}^2-a_{16}a_{34}^2-a_{16}a_{35}^2+a_{16}a_{36}^2),
\end{array}
$$
$$
\begin{array}{rcl}
H_{2,4}  & = & \frac{-c^2}{\sqrt{-\lambda}}(a_{14}^2a_{24}+2a_{14}a_{15}a_{25}+2a_{14}a_{16}a_{26}-a_{15}^2a_{24}-a_{16}^2a_{24}+a_{24}^3+a_{24}a_{25}^2+a_{24}a_{26}^2\\
               & &+a_{24}a_{34}^2-a_{24}a_{35}^2-a_{24}a_{36}^2+2a_{25}a_{34}a_{35}+2a_{26}a_{34}a_{36}),\\
H_{2,5}  & = & \frac{c^2}{\sqrt{-\lambda}}(a_{14}^2a_{25}-2a_{14}a_{15}a_{24}-a_{15}^2a_{25}-2a_{15}a_{16}a_{26}+a_{16}^2a_{25}-a_{24}^2a_{25}-2a_{24}a_{34}a_{35}\\
               & &-a_{25}^3-a_{25}a_{26}^2+a_{25}a_{34}^2-a_{25}a_{35}^2+a_{25}a_{36}^2-2a_{26}a_{35}a_{36}),\\
H_{2,6}  & = & \frac{c^2}{\sqrt{-\lambda}}(a_{14}^2a_{26}-2a_{14}a_{16}a_{24}+a_{15}^2a_{26}-2a_{15}a_{16}a_{25}-a_{16}^2a_{26}-a_{24}^2a_{26}-2a_{24}a_{34}a_{36}\\
               & &-a_{25}^2a_{26}-2a_{25}a_{35}a_{36}-a_{26}^3+a_{26}a_{34}^2+a_{26}a_{35}^2-a_{26}a_{36}^2),\\
H_{3,4}  & = & \frac{-c^2}{\sqrt{-\lambda}}(a_{14}^2a_{34}+2a_{14}a_{15}a_{35}+2a_{14}a_{16}a_{36}-a_{15}^2a_{34}-a_{16}^2a_{34}+a_{24}^2a_{34}+2a_{24}a_{25}a_{35}\\
              & &+2a_{24}a_{26}a_{36}-a_{25}^2a_{34}-a_{26}^2a_{34}+a_{34}^3+a_{34}a_{35}^2+a_{34}a_{36}^2),\\
H_{3,5}  & = & \frac{c^2}{\sqrt{-\lambda}}(a_{14}^2a_{35}-2a_{14}a_{15}a_{34}-a_{15}^2a_{35}-2a_{15}a_{16}a_{36}+a_{16}^2a_{35}+a_{24}^2a_{35}-2a_{24}a_{25}a_{34}\\
               & &-a_{25}^2a_{35}-2a_{25}a_{26}a_{36}+a_{26}^2a_{35}-a_{34}^2a_{35}-a_{35}^3-a_{35}a_{36}^2),\\
H_{3,6}  & = & \frac{c^2}{\sqrt{-\lambda}}(a_{14}^2a_{36}-2a_{14}a_{16}a_{34}+a_{15}^2a_{36}-2a_{15}a_{16}a_{35}-a_{16}^2a_{36}+a_{24}^2a_{36}-2a_{24}a_{26}a_{34}\\
               & &+a_{25}^2a_{36}-2a_{25}a_{26}a_{35}-a_{26}^2a_{36}-a_{34}^2a_{36}-a_{35}^2a_{36}-a_{36}^3),
\end{array}
$$
where $H_{i,j} = h(e_i,e_j)$.
Observe that up to multiplication by $\sqrt{-\lambda}$, the nonzero terms are all homogeneous polynomials of third degree in the $a_{ij}$.

We are looking for coupled structures inducing either the standard metric or the Jensen metric which with respect to the considered basis can be written as the identity matrix and as \eqref{Jensen}, respectively. Thus, since $\omega\W\psip=0, d\psip = 0$ and $d\omega^2=0$, we first have to solve the system obtained by imposing that the matrix $H$ is proportional to the identity matrix or to the matrix \eqref{Jensen} under the assumption $\lambda<0$ and then, if we find solutions of this system, we need to impose that the normalization condition is satisfied in order to obtain what we want.\\
\noindent
{\bf Case 1: the standard metric}\\
Since rescaling a metric with a positive coefficient does not change the Ricci tensor, we are looking for solutions of the equation
$$H = \alpha I,$$
where $\alpha$ is a positive real number. 

Since the entries in the diagonal of $H$ are all equal, we only have to solve the system of equations
$$H_{i,j} = 0, i=1,2,3, j = 4,5,6$$
under the assumptions $H_{1,1}\neq 0$ and $\lambda<0$.

For $i,j = 1,\ldots,6$, we let
$$\tilde{H}_{i,j} := \sqrt{-\lambda}H_{i,j}.$$
Then, as already observed, the $\tilde{H}_{i,j}$ are homogeneous polynomials of degree 3 in the $a_{ij}$ and under our assumptions $H_{i,j} = 0$ 
if and only if $\tilde{H}_{i,j} = 0$ for $ i=1,2,3, j = 4,5,6$.

Since we have a system of equations involving homogeneous polynomials of the same degree and we are looking for solutions defined up to a multiplicative constant, 
let us consider the projective space $\mathbb{CP}^8$ with coordinate ring $$\C[a_{14},a_{15},a_{16},a_{24},a_{25},a_{26},a_{34},a_{35},a_{36}]$$
and the homogeneous ideals
$$P:= \langle\tilde{H}_{1,1}\rangle,$$
$$Q:=\langle\tilde{H}_{1,4},\tilde{H}_{2,4},\tilde{H}_{3,4},\tilde{H}_{1,5},\tilde{H}_{2,5},\tilde{H}_{3,5},\tilde{H}_{1,6},\tilde{H}_{2,6},\tilde{H}_{3,6}\rangle.$$
What we are looking for is the set of points $[a_{14}:\ldots:a_{36}]$ lying in the projective variety $V(Q)$ but not in $V(P)$ and for which $\lambda<0$. 
It is known that 
$$\overline{V(Q)-V(P)} \subseteq V(Q:P),$$
where $Q:P$ is the ideal quotient of $Q$ by $P$. 
In our case it is possible to show that $Q:P = \langle a_{14},a_{15},a_{16},a_{24},a_{25},a_{26},a_{34},a_{35},a_{36}\rangle,$ therefore $V(Q:P) =  \emptyset$.

We have then proved that on $S^3\times S^3$ there are no invariant coupled structures inducing the standard metric.\\
\noindent
{\bf Case 2: the Jensen metric}\\
Following the same idea of the previous case and looking at the entries of the matrix \eqref{Jensen}, we have now to consider the ideals $P$ and 
$$R:= \langle\tilde{H}_{1,5},\tilde{H}_{1,6},\tilde{H}_{2,4},\tilde{H}_{2,6},\tilde{H}_{3,4},\tilde{H}_{3,5},\tilde{H}_{2,5}-\tilde{H}_{3,6},\tilde{H}_{3,6}-\tilde{H}_{1,4},\tilde{H}_{1,1} +2\tilde{H}_{1,4}\rangle$$
and look for those points lying in the projective variety $V(R)$ but not in $V(P)$ and for which $\lambda<0$.
Now 
$$R:P = \langle a_{15}, a_{16},a_{24},a_{26},a_{34},a_{35},a_{25}-a_{14}, a_{36}-a_{14}\rangle,$$
then 
$$V(R:P) = \{[\gamma: 0: 0: 0: \gamma: 0: 0: 0: \gamma] : \gamma\in\C-\{0\} \}$$
is a point in $\mathbb{CP}^8$ and since $\C$ is algebrically closed and $R$ is a radical ideal
$$V(R:P) = \overline{V(R)-V(P)}.$$
Moreover the requested condition on $\lambda$ is satisfied, indeed:
$$
\lambda  =  -3c^4\gamma^4 < 0.
$$
The coupled structures we are interested in are obtained when $\gamma$ is a negative real number, in this case we have:
\begin{eqnarray}
\omega & = & \gamma (e^{14} +  e^{25} +  e^{36}),\nonumber \\
\psip & = & c\gamma (e^{234}- e^{156}- e^{135}+ e^{246}+ e^{126}- e^{345}),\nonumber \\
\psim &=&  \frac{c\gamma}{\sqrt{3}}(2 e^{123}-e^{126}+e^{135}-e^{156}-e^{234}+e^{246}-e^{345}+2e^{456}).\nonumber
\end{eqnarray}
The forms $\omega$ and $\psip$ are stable and the normalization condition implies 
$$c = \pm \sqrt{\frac{-2\gamma}{\sqrt{3}}},$$
in both cases we have a nearly K\"ahler structure.
\endproof

Anyway, it is not difficult to find examples of left invariant coupled structures on $S^3\times S^3$:
\begin{ex}
The pair
$$\renewcommand\arraystretch{1.4}
\begin{array}{lcl}
\omega &=& -\sqrt{3}e^{16}-e^{24}-e^{25}-e^{35},\\ 
\psip &=& \sqrt[4]{3}(-\sqrt{3}e^{236}+\sqrt{3}e^{145}+e^{134}+e^{256}+e^{135}-e^{246}-e^{125}-e^{346}),
\end{array}
$$
defines a coupled structure on $\mathfrak{su}(2)\oplus\mathfrak{su}(2)$ such that $\psip = \sqrt[4]{3}d\omega$. In particular, we have that  $(h,J,\omega)$, where $J=J_{\psi_+}$, 
defines a quasi K\"ahler structure on $\mathfrak{su}(2)\oplus\mathfrak{su}(2)$.
\end{ex}

\section{Einstein half-flat structures on 6-solvmanifolds}\label{secsol}
In section \ref{s3s3} we saw that on the compact manifold  $S^3\times S^3$ there are examples of left invariant half-flat, double and nearly K\"ahler structures inducing an Einstein metric and that  there are no left invariant coupled structures inducing any of the Einstein metrics known up to now. In this section we turn our attention to the noncompact homogeneous case.

It is well known that noncompact homogeneous Einstein manifolds have non positive Ricci curvature. Moreover, if the Ricci curvature is zero, then the Riemannian metric is flat (\cite{AK}) and the manifold is isometric to the product of a flat torus and the Euclidean space. 
Up to now the only known examples of noncompact homogeneous Einstein manifolds are Einstein solvmanifolds, that is simply connected solvable Lie groups endowed with a left invariant Einstein metric. It has been conjectured by Alekseevskii that any noncompact homogeneous Einstein manifold might be of this kind (see \cite{Bes} for the statement of the conjecture and refer to the recent works \cite{AL} and \cite{JP} and the references therein for more details and the most recent results on it).

Let $(S,h)$ be a solvmanifold, we can identify the left invariant metric $h$ on the simply connected solvable Lie group $S$ with the inner product $h_0$ determined by it on the solvable Lie algebra $\mathfrak{s}$ of $S,$ the pair $(\mathfrak{s},h_0)$ is said a {\bf metric solvable Lie algebra}. Two metric Lie algebras $(\mathfrak{s},h_0)$ and $(\mathfrak{s}',h'_0)$ are {\bf isomorphic} if there exists a Lie algebra isomorphism $F:\mathfrak{s}\rightarrow\mathfrak{s}'$ which is also an isometry of Euclidean spaces.

In \cite{Lau}, Lauret showed that any Einstein solvmanifold is {\bf standard}, i.e.~the orthogonal complement $\mathfrak{a}$ to $\mathfrak{n} = [\mathfrak{s},\mathfrak{s}]$ is always an abelian subalgebra of $\mathfrak{s} = \mathfrak{n}\oplus\mathfrak{a}$. The dimension of $\mathfrak{a}$ is the {\bf algebraic rank} of $\mathfrak{s}$. 

The properties of standard Einstein solvmanifolds were studied by Heber in \cite{He}. In particular, he showed that a standard Einstein metric is unique up to isometries and scaling among left invariant metrics and that, up to metric Lie algebra isomorphisms, a standard metric Lie algebra with Einstein inner product is an Iwasawa-type algebra.

It was proved by Dotti in \cite{Dot} that solvmanifolds $(S,h)$ with unimodular solvable Lie group $S$ and Einstein metric $h$ are flat.

The 6-dimensional Einstein solvmanifolds were classified by Nikitenko and Nikonorov in \cite{NN}. The result is recalled in the next Theorem. 
Instead of the Lie algebra structure equations given (in the original formulation of the Theorem) by the nontrivial Lie brackets of the basis vectors, we write here the structure equations in terms of the Chevalley-Eilenberg differential of the basis 1-forms, since we will use these in our next computations.

\begin{teo}[\cite{NN}]\label{NNthm}
Let $(\mathfrak{s}, h)$ be a 6-dimensional nonunimodular metric solvable Lie algebra with Einstein inner product $h$ such that ${\rm Ric}(h)=-r^2h$, 
where $r>0$. Then $(\mathfrak{s}, h)$ is isomorphic to one of the metric Lie algebras contained in Table 1. For each algebra, $(e_1,\ldots,e_6)$ is an h-orthonormal basis with dual basis $(e^1,\ldots,e^6)$. 
\end{teo}

\begin{table}[h]
\caption{6-dimensional nonunimodular metric solvable Lie algebras with Einstein inner product $h = \sum_{i=1}^6(e^i)^2$. The Lie algebra $\mathfrak{s}_{10}$ depends on a parameter 
$0\leq t \leq \frac{1}{\sqrt{22}}$ and $\mathfrak{s}_{12}$ depends on two parameters $0\leq s \leq t \leq 1$.} 
\centering
\renewcommand\arraystretch{1.8}
\begin{tabular}{|c|c|}
\hline
$\mathfrak{s}_\cdot$&{\rm Structure equations} $(de^1,de^2,de^3,de^4,de^5,de^6)$ \\ \hline
$\mathfrak{s}_1$      & ${\scriptstyle \left(\frac{r}{2\sqrt{2}}e^{16},\frac{r}{2\sqrt{2}}e^{26},\frac{r}{2\sqrt{2}}e^{36},\frac{r}{2\sqrt{2}}e^{46},-\frac{r}{\sqrt{2}}e^{12}-\frac{r}{\sqrt{2}}e^{34}+\frac{r}{\sqrt{2}}e^{56},0\right)}$         \\ \hline 
$\mathfrak{s}_2$       & ${\scriptstyle \left(2r\sqrt{\frac{2}{105}}e^{16},r\sqrt{\frac{3}{70}}e^{26},-\frac{2r}{\sqrt{7}}e^{12}+r\sqrt{\frac{7}{30}}e^{36},2r\sqrt{\frac{3}{70}}e^{46},-r\sqrt{\frac27}e^{14}-\frac{2r}{\sqrt{7}}e^{23}+r\sqrt{\frac{10}{21}}e^{56}, 0\right)}$           \\ \hline		
$\mathfrak{s}_3$       & ${\scriptstyle \left(\frac{r}{\sqrt{55}}e^{16},\frac{2r}{\sqrt{55}}e^{26},-r\sqrt{\frac{6}{11}}e^{12}+\frac{3r}{\sqrt{55}}e^{36},-r\sqrt{\frac{6}{11}}e^{13}+\frac{4r}{\sqrt{55}}e^{46},-\frac{2r}{\sqrt{11}}e^{14} -\frac{2r}{\sqrt{11}}e^{23}+\frac{5r}{\sqrt{55}}e^{56},0\right)}$                                           \\ \hline
$\mathfrak{s}_4$       &${\scriptstyle \left(\frac{r\sqrt{6}}{30}e^{16},\frac{3r\sqrt{6}}{20}e^{26},-\frac{r}{\sqrt{2}}e^{12}+\frac{11r\sqrt{6}}{60}e^{36},-r\sqrt{\frac23}e^{13}+\frac{13r\sqrt{6}}{60}e^{46},-\frac{r}{\sqrt{2}}e^{14}+\frac{r\sqrt{6}}{4}e^{56},0\right)}$                    \\ \hline
$\mathfrak{s}_5$       & ${\scriptstyle \left( \frac{r}{3\sqrt{2}}e^{16},\frac{r}{2\sqrt{2}}e^{26},\frac{r}{2\sqrt{2}}e^{36},-\frac{r}{\sqrt{2}}e^{12}+ \frac{5r}{6\sqrt{2}}e^{46},-\frac{r}{\sqrt{2}}e^{13}+ \frac{5r}{6\sqrt{2}}e^{56},0 \right)}$                             \\ \hline
$\mathfrak{s}_{6}$    & ${\scriptstyle \left(\frac{r}{2\sqrt{6}}e^{16},\frac{r}{2\sqrt{6}}e^{26},-r\sqrt{\frac23}e^{12}+\frac{r}{\sqrt{6}}e^{36},-\frac{r}{\sqrt{2}}e^{13}+r\frac{\sqrt{6}}{4}e^{46},-\frac{r}{\sqrt{2}}e^{23}+r\frac{\sqrt{6}}{4}e^{56},0\right)}$    			             \\ \hline
$\mathfrak{s}_{7}$    & ${\scriptstyle \left(\frac{r}{\sqrt{39}}e^{16},\frac{2r}{\sqrt{39}}e^{26},-r\sqrt{\frac23}e^{12} + \frac{3r}{\sqrt{39}}e^{36},-r\sqrt{\frac23}e^{13} +\frac{4r}{\sqrt{39}}e^{46},\frac{3r}{\sqrt{39}}e^{56},0\right)}$    	                     \\ \hline
$\mathfrak{s}_{8}$    & ${\scriptstyle \left(r\sqrt{\frac{2}{21}}e^{16},r\sqrt{\frac{2}{21}}e^{26},-r\sqrt{\frac23}e^{12}+ 2r\sqrt{\frac{2}{21}}e^{36},r\sqrt{\frac{3}{14}}e^{46},r\sqrt{\frac{3}{14}}e^{56},0\right)}$   \\ \hline
$\mathfrak{s}_{9}$    &${\scriptstyle \left(\frac{r}{\sqrt{5}}e^{16},\frac{r}{\sqrt{5}}e^{26},\frac{r}{\sqrt{5}}e^{36},\frac{r}{\sqrt{5}}e^{46},\frac{r}{\sqrt{5}}e^{56},0\right)}$   \\ \hline
$\mathfrak{s}_{10}$    & ${\scriptstyle \left(\frac{2r}{\sqrt{33}}e^{15}+rt e^{16}+r\sqrt{\frac12-11t^2} e^{26},\frac{2r}{\sqrt{33}}e^{25}+r\sqrt{\frac12-11t^2} e^{16} + rt e^{26},
-r\sqrt{\frac23}e^{12}+\frac{4r}{\sqrt{33}}e^{35}+2rt e^{36},\frac{3r}{\sqrt{33}}e^{45}-4rt e^{46},0,0\right)}$ 							     		                          \\ \hline
$\mathfrak{s}_{11}$    & ${\scriptstyle \left(\frac{r}{\sqrt{30}}e^{15} + \frac{3r}{\sqrt{30}}e^{16},\frac{2r}{\sqrt{30}}e^{25} - \frac{4r}{\sqrt{30}}e^{26},-r\sqrt{\frac23}e^{12}+\frac{3r}{\sqrt{30}}e^{35} 
- \frac{r}{\sqrt{30}}e^{36},-r\sqrt{\frac23}e^{13}+\frac{4r}{\sqrt{30}}e^{45} + \frac{2r}{\sqrt{30}}e^{46},0,0\right)}$ 							     	                                          \\ \hline
$\mathfrak{s}_{12}$    & ${\scriptstyle \left(\frac{r}{2}e^{15} + r\frac{1+s+t}{2\sqrt{1+t^2+s^2}}e^{16},\frac{r}{2}e^{25} + r\frac{1-s-t}{2\sqrt{1+t^2+s^2}}e^{26},
\frac{r}{2}e^{35} + r\frac{t-s-1}{2\sqrt{1+t^2+s^2}}e^{36},\frac{r}{2}e^{45} + r\frac{s-t-1}{2\sqrt{1+t^2+s^2}}e^{46},0,0\right)}$ 			\\ \hline
$\mathfrak{s}_{13}$  & ${\scriptstyle\left(\frac{r}{\sqrt{3}}e^{14} - \frac{2r}{\sqrt{6}}e^{16},\frac{r}{\sqrt{3}}e^{24} + \frac{r}{\sqrt{2}}e^{25} + \frac{r}{\sqrt{6}}e^{26},\frac{r}{\sqrt{3}}e^{34} - \frac{r}{\sqrt{2}}e^{35} + \frac{r}{\sqrt{6}}e^{36}, 0,0,0\right)}$                  \\ \hline		                
\end{tabular}
\end{table}
\renewcommand\arraystretch{1}

\begin{remark}
All the solvable metric Lie algebras appearing in Table 1 are of Iwasawa-type. It is worth emphasizing here that for each of the $\mathfrak{s}_i$, the inner product, with respect to which the basis $(e_1,\ldots,e_6)$ is orthonormal, is the only one having the property of being Einstein. This follows from the result of Heber and the fact that two solvable metric Lie algebras of Iwasawa-type are isomorphic if and only if the corresponding solvmanifolds are isometric as Riemannian manifolds. 
\end{remark}

As a consequence of the previous Theorem, for the 6-dimensional homogeneous Einstein manifolds of nonpositive sectional curvature they showed what follows.
\begin{teo}[\cite{NN}]\label{NNthm2}
Let $(M,h)$ be a 6-dimensional, connected, simply connected, homogeneous Einstein manifold of nonpositive sectional curvature, then it is symmetric or isometric to one of the solvmanifolds 
of negative sectional curvature generated by the metric Lie algebras $\mathfrak{s}_5$, $\mathfrak{s}_8$. Moreover, in the symmetric case $(M,h)$ is obtained as the solvmanifold corresponding to 
the metric Lie algebras $\mathfrak{s}_1, \mathfrak{s}_9, \mathfrak{s}_{10}$ for $t=\frac{1}{\sqrt{22}}, \mathfrak{s}_{11}, \mathfrak{s}_{12}$ for $(s,t) = (0,0)$ and $(s,t)=(1,1)$ and $\mathfrak{s}_{13}$.   
\end{teo}

Now we focus on the problem of finding left invariant half-flat structures on 6-dimensional Einstein solvmanifolds inducing the Einstein (non Ricci-flat) metric, these are in 1-1 correspondence with half-flat structures inducing the Einstein inner product on 6-dimensional nonunimodular solvable metric Lie algebras. 

It is worth recalling here an obstruction to the existence of half-flat structures on 6-dimensional Lie algebras shown by Freibert and Schulte-Hengesbach in \cite{FSD}. 
This result is a refinement of the one obtained by Conti in \cite{Con}.

\begin{prop}[\cite{FSD}]\label{obst}
Let $\mathfrak{g}$ be a 6-dimensional Lie algebra with volume form $\Omega\in\Lambda^6(\mathfrak{g}^*)$. If there exists a nonzero $\alpha\in\mathfrak{g}^*$ such that
$$\alpha\W\tilde{J}_\rho^*\alpha\W\sigma = 0$$
for all closed 3-forms $\rho\in\Lambda^3(\mathfrak{g}^*)$ and closed 4-forms $\sigma\in\Lambda^4(\mathfrak{g}^*)$, where for any $X\in\mathfrak{g}$
$$\tilde{J}_\rho^*\alpha(X)\Omega = \alpha\W(i_X\rho)\W\rho,$$
then $\mathfrak{g}$ does not admit any half-flat $SU(3)$-structure.
\end{prop}

In \cite{FSD} and \cite{FSI} the authors completely classified the left invariant half-flat structures on 6-dimensional decomposable Lie groups (using also the classification contained in \cite{SPT}) and on 6-dimensional indecomposable Lie groups with 5-dimensional nilradical. These classifications will be useful in the proof of the following 
\begin{teo}\label{mainthm}
There are no half-flat $\SU(3)$-structures inducing the Einstein metric on the rank 1 solvable metric Lie algebras $\mathfrak{s}_i, i=1\ldots9$, and on the rank 2 solvable metric Lie algebra $ \mathfrak{s}_{12}$
and there are no coupled $\SU(3)$-structures inducing the Einstein metric on the rank 2 metric Lie algebras $\mathfrak{s}_{10}, \mathfrak{s}_{11}$ and on the rank 3 metric Lie algebra $\mathfrak{s}_{13}$.
\end{teo}
\proof We will prove the theorem as follows: in the list of Einstein solvable metric Lie algebras we first exclude the ones that do not admit a half-flat structure using the results of \cite{FSD} and \cite{FSI}, then we will show the result by direct computations in the remaining cases.

The rank 1 Lie algebra $\mathfrak{s}_9$ is indecomposable and has abelian nilradical, therefore it does not admit any half-flat structure by Proposition 4 of \cite{FSI}.
By Theorem 2 of \cite{FSI} we have that the Lie algebras $\mathfrak{s}_i$ with $i = 1,2,4,5,7,8$ do not admit any half-flat structure since they are isomorphic to the Lie algebras 
$A^{1,0,0}_{6,82}, A^{{4/3}}_{6,94}, A^{{9/2}}_{6,71}, A^{{2/5,1}}_{6,54}, A^{3,2}_{6,39}, A^{{2/3,2/3,1}}_{6,13}$ respectively, whereas the Lie algebras $\mathfrak{s}_3,\mathfrak{s}_6$ 
admit a half-flat structure since the former is isomorphic to $A_{6,99}$ and the latter to $A^1_{6,76}$. By Theorem 1 of \cite{FSD} we also have that on $\mathfrak{s}_{13}$ there exist half-flat structures 
since it is isomorphic to $\mathfrak{r}_2\oplus\mathfrak{r}_2\oplus\mathfrak{r}_2$.
Finally, applying Proposition \ref{obst} to $\mathfrak{s}_{12}$ with $\alpha = e^6$ we obtain that this 2-parameter family of Lie algebras does not admit any half-flat structure.

We can now start with the second part of the proof. For $k=3,6,$ let $\omega\in\Lambda^2(\mathfrak{s}^*_k)$ and $\psip\in\Lambda^3(\mathfrak{s}^*_k)$ be generic forms, with respect to the basis $(e^1,\ldots,e^6)$ given in Theorem \ref{NNthm} we can write
\begin{equation}\label{omggen}
\begin{array}{rcl}
\omega & = & b_1e^{12} + b_2e^{13}+b_3e^{14}+b_4e^{15} + b_5e^{16} + b_6e^{23} + b_7 e^{24} +b_8e^{25}\\
               &  &  +b_9e^{26} + b_{10}e^{34} + b_{11}e^{35} + b_{12}e^{36} +b_{13}e^{45} + b_{14} e^{46} + b_{15}e^{56}
\end{array}
\end{equation}
and 
\begin{equation}\label{psipgen}
\begin{array}{rcl}
\psip & = & a_{1}e^{123}+a_{2}e^{124}+a_{3}e^{125}+a_{4}e^{126}+a_{5}e^{134}+a_{6}e^{135}+a_{7}e^{136}\\
         & &+a_{8}e^{145}+a_{9}e^{146}+a_{10}e^{156}+a_{11}e^{234}+a_{12}e^{235}+a_{13}e^{236}+a_{14}e^{245}\\
          & &+a_{15}e^{246}+a_{16}e^{256}+a_{17}e^{345}+a_{18}e^{346}+a_{19}e^{356}+a_{20}e^{456},
\end{array}
\end{equation}
where $a_i$ and $b_j$ are real constants. 
Moreover, we denote by $\beta_{i_1\ldots i_5}$ and $\gamma_{i_1\ldots i_5}$ the components of the 5-forms $\omega\W\psip$ and $d\omega^2$ respectively, so that
$$\renewcommand\arraystretch{1.8}
\begin{array}{rcl}
\omega\W\psip &=& \sum_{1\leq i_1<i_2<\ldots<i_5\leq6}\beta_{i_1\ldots i_5}e^{i_1\ldots i_5},\\
d\omega^2 &=& \sum_{1\leq i_1<i_2<\ldots<i_5\leq6}\gamma_{i_1\ldots i_5}e^{i_1\ldots i_5}.
\end{array}
$$
Observe that the non vanishing $\beta$ are always homogeneous polynomials of degree 2 in the $a_i,b_j$, while the non vanishing $\gamma$ are always homogeneous polynomials of degree 2 in the $b_j$.

For each Lie algebra $\mathfrak{s}_k$ we impose the conditions the forms \eqref{omggen} and \eqref{psipgen} have to satisfy in order to be a half-flat $\SU(3)$-structure inducing the Einstein metric. What we have to do is to solve the equations obtained from 
\begin{equation}\label{systemSU3}
\begin{cases}
\omega\W\psip = 0\\
d\psip = 0\\
d\omega^2 = 0
\end{cases}
\end{equation}
under the assumptions $\lambda=\lambda(\psip)<0, \omega^3\neq0$. Moreover, since we are considering a basis which is orthonormal with respect to the Einstein metric (see Theorem \ref{NNthm}), 
we have also to impose that the entries $H_{i,j}=h(e_i,e_j)$ of the matrix $H$ associated to $h(\cdot,\cdot) = \omega(\cdot, J_{\psip}\cdot)$ with respect to the basis $(e_1,\ldots,e_6)$ satisfy
\begin{equation}\label{hidentity}
\begin{array}{rcl}
H_{i,j} &=&  0,\ \mbox{for\ } 1\leq i,j \leq 6 \ \mbox{and\ } i\neq j,\\
H_{i,i}-H_{i+1,i+1} &=& 0,\ \mbox{for\ } 1\leq i \leq5,
\end{array}
\end{equation}
with $H_{i,i}>0$. These conditions give us a system of polynomial equations in 35 unknowns to solve under some constraints on them we will specify case by case. 
Since the expressions of the unknowns we obtain solving the equations are often too long to be written down, in what follows we will point out only from which equation a certain unknown is obtained, specifying its value only if it is zero.

Let us start with the Lie algebra $\mathfrak{s}_6$, whose structure equations are given in Table 1. 
We solve all the linear equations in the $a_i$ deriving from $d\psip=0$. Then looking at the expression of $\lambda$ we deduce that $a_6\neq0$.
We can then solve all the equations obtained from $\omega\W\psip=0$, $d\omega^2=0$ and $H_{i,j} = 0$ for $i\neq j$ using $a_6\neq0$ and comparing case by case each equation with $H_{i,i}$ and $\omega^3$.
After doing this, $H$ becomes a diagonal matrix and we have to solve the remaining 5 equations of \eqref{hidentity}, which do not have any solution under the constraints $H_{i,i}\neq0$ and $\lambda\neq0$. 

For the Lie algebra $\mathfrak{s}_3$ we can argue in a similar way, but instead of working on it, we can show the result on the Lie algebra $A_{6,99}\cong\mathfrak{s}_3$ since the computations are less involved. The structure equations of $A_{6,99}$ are given for example in \cite{FSI}, with respect to a basis $(\ff_1,\ldots,\ff_6)$ with dual basis $(\ff^1,\ldots,\ff^6)$ they are
$$\left(5\ff^{16}+\ff^{25}+\ff^{34},4\ff^{26}+\ff^{35},3\ff^{36}+\ff^{45},2\ff^{46},\ff^{56},0\right).$$
We consider the generic forms $\omega$ and $\psip$ as in \eqref{omggen} and \eqref{psipgen} with $e^i$ replaced by $\ff^i$. 
Observe that the matrix $H$ associated to the Einstein inner product with respect to the basis $(\ff_1,\ldots,\ff_6)$ is not proportional to the identity but it is still diagonal, thus we still have to solve the equations $H_{i,j}=0$ for $i\neq j$. First of all we solve the linear equations in the $a_i$ obtained from $d\psip=0$. 
Then we observe that having $b_1=0$ or $a_6=0$ leads to a contradiction after solving some equations: if $b_1=0$ we can use $H_{1,1},H_{2,2}\neq0$ to solve $H_{1,2}=0$, $H_{1,3}=0$, $\beta_{12345}=0$, $\gamma_{12346}=0$, $\gamma_{12356}=0$, $\beta_{12346}=0$, but then $\gamma_{12456}$ can not be zero; if $a_6=0$ we use $H_{1,1},H_{2,2}\neq0$ to solve $H_{1,2}=0$, $H_{1,3}=0$, $\beta_{12345}=0$, $H_{1,5}=0$, $\beta_{12356}=0$, $\gamma_{12356}=0$, $H_{2,3}=0$, $\beta_{12346}=0$ and obtain that $\gamma_{12346}=0$ if and only if $H_{1,1}H_{3,3}=0$.
Thus we assume  $b_1\neq0$ and $a_6\neq0$. Under these constraints and comparing case by case the polynomial we want to be zero with $H_{i,i}$ and $\lambda$, we can get the expression of $b_4$ from the equation $H_{1,2}=0$, $b_9$ from $H_{2,3}=0$, $b_7$ from $\beta_{12345}=0$, $b_{10}$ from $\gamma_{12346}=0$, 
$b_{11}$ from $\gamma_{12356}=0$, $a_{18}$ from $\beta_{12346}=0$, $b_2=0$ from $\beta_{12356}=0$, $b_6=0$ from $H_{1,3}=0$, $a_8=0$ from $H_{3,4}=0$, $a_{14}=0$ from $H_{3,5}=0$, $a_{17}=0$ from $H_{3,6}=0$, $a_{10}=0$ from $H_{1,4}=0$, $a_{19}=0$ from $H_{1,6}=0$, $b_{14}=0$ from $\beta_{13456}=0$, $a_{20}$ from $H_{1,5}=0$, $b_3$ from $H_{2,4}=0$, $b_8$ from $H_{2,6}=0$. Now $H_{4,6}=0$ implies $\omega^3=0$.
 
We can now turn our attention to the Lie algebras $\mathfrak{s}_{10}, \mathfrak{s}_{11}$ and $\mathfrak{s}_{13}$, we will show that none of these admits a coupled structure inducing the Einstein metric. The way in which we proceed is similar to the one followed for $\mathfrak{s}_6$ and $\mathfrak{s}_3$, but in this case we will consider a generic $\omega$ 
of the form \eqref{omggen} and $\psip = cd\omega$ for $c\in\R-\{0\}$. Observe that the second condition of \eqref{systemSU3} is satisfied since $\psip$ is now an exact 3-form and that the first and the third condition are actually the same. For each Lie algebra we consider the structure equations given in Table 1.
 
Consider $\mathfrak{s}_{10}$, this is a 1-parameter family of Lie algebras depending on $t\in\left[0,\frac{1}{\sqrt{22}}\right]$. Since $H_{3,3}$ can not be zero, we have that $b_{10}\neq0$, $b_2\neq\pm b_6$ and $t\neq\frac{1}{\sqrt{22}}$. The way in which we solve the equations depends on whether $t=\frac{7}{2\sqrt{330}}$ or not. If $t\neq\frac{7}{2\sqrt{330}}$ we can use $b_{10}\neq0$ to obtain $b_1$ from $\gamma_{12345}=0$, $b_{15}$ from $\gamma_{12456}=0$, $b_{5}$ from $\gamma_{13456}=0$ and $b_{9}$ from $\gamma_{23456}=0$.  
Then $b_{12}$ from $H_{3,4}=0$, $b_3=0$ from $H_{1,3}=0$, $b_7=0$ from $H_{2,3}=0$, $b_{11}=0$ from $H_{4,5}=0$, $b_{8}$ from $H_{1,4}=0$, $b_{4}$ from $H_{2,4}=0$,  $b_{13}=0$ from $H_{3,5}=0$. Now $H_{3,6}=0$ if and only if $\lambda=0$. If $t=\frac{7}{2\sqrt{330}}$ the computations are the same until we arrive to the equation $H_{2,4}=0,$ 
which has no solutions since $H_{2,4}$ is proportional to $\lambda$.

For $\mathfrak{s}_{11}$ we have that $b_{10}\neq0$ and $\sqrt{5}b_{3}b_{7}-b_{10}b_{13}-3b_{10}b_{14}\neq0$, otherwise $H_{4,4}=0$. 
Using $b_{10}\neq 0$ we can obtain $b_1,b_4,b_8,b_{15}$ from $\gamma_{12345}=0,\gamma_{13456}=0,\gamma_{23456}=0,\gamma_{12456}=0$, respectively, $b_{11}$ from $H_{3,4}=0$, $b_{9}$ from $H_{2,4}=0$ and $b_{5}$ from $H_{1,4}=0$. 
Then using also the other constraint we get $b_{12}$ from $H_{4,6}=0$. Now $H_{4,5}=0$ if and only if $b_7=0$ or $b_{14}=-2b_{13}$. 
If $b_7=0$, from $H_{2,5}=0$ and $\lambda\neq0$ we have $b_3=0$ but then $H_{1,2}=0$ only if either $\lambda=0$ or $H_{1,1}=0$. Thus $b_7\neq0$ and $b_{14}=-2b_{13}$. 
Moreover, $b_6\neq0$, otherwise $\lambda$ would be proportional to $H_{2,3}$. Thus we can solve $H_{2,3}=0$ to get $b_3$ and use $\lambda\neq0$ to solve $H_{1,2}=0$ and obtain $b_6$. Now $H_{3,5}$ is proportional to $\lambda$, therefore it can not be zero. 

In the last case $\mathfrak{s}_{13}$ we can see that $b_1, b_2, b_6\neq0$ and $b_3\neq \sqrt{2}{}b_5$ from the fact that the entries in the diagonal of $H$ can not be zero. 
Solving the equations $\gamma_{12456}=0,\gamma_{12345}=0,\gamma_{12346}=0,\gamma_{13456}=0,\gamma_{23456}=0$ under the previous constraints we obtain the expressions of $b_{14},b_8,b_9,b_{13},b_{15},$ respectively.
Then we get $b_{12}$ from $H_{1,2}=0$, $b_{10}$ from $H_{1,3}=0$, $b_{11}$ from $H_{1,5}=0$, $b_5$ from $H_{2,5}=0$ and $b_4=0$ from $H_{2,6}=0$. Now $H_{2,4}=0$ if and only if $\lambda=0$.
\endproof

\begin{remark}
In this case it is in principle possible to use the properties of algebraic varieties to find solutions as we did in the proof of Theorem \ref{nocpds3}. However the computations here are more involved since we have more unknowns (35 or 15 instead of 9) and more equations arising from the fact that some defining conditions for an $\SU(3)$-structure that were easily verified in the case of $S^3\times S^3$ have to be imposed in this case. 
\end{remark}

From the fact that the class of coupled structures is a subclass of the half-flat one, we can use the result of the previous theorem together with Theorem \ref{NNthm} to obtain: 
\begin{corol}
Let $(\mathfrak{s}, h)$ be a 6-dimensional nonunimodular solvable metric Lie algebra with $h$ Einstein. Then on $\mathfrak{s}$ there are no coupled $\SU(3)$-structures inducing the Einstein inner product. 
\end{corol}

Moreover, from the previous theorem and the Theorem \ref{NNthm2} we obtain a constraint for the existence of coupled structures inducing Einstein metrics on homogeneous spaces:
\begin{corol}
Let $(M,h)$ be a 6-dimensional, connected, simply connected, homogeneous Einstein manifold of nonpositive sectional curvature. Then there are no left invariant coupled $\SU(3)$-structures on $M$ inducing the Einstein metric.
\end{corol}

\end{document}